\documentclass[12pt,a4paper]{article}

\usepackage{amsfonts}
\usepackage{amssymb}
\usepackage{amsthm}
\usepackage{amsmath}
\usepackage{enumerate}
\usepackage{hyperref}

\usepackage{graphicx}

\def\be{\begin{equation}}
\def\ee{\end{equation}}
\def\IC{\mathbb{C}}
\def\ID{\mathbb{D}}
\def\CD{\overline{\mathbb{D}}}

\def\IZ{\mathbb{Z}}

\def\IT{\mathbb{T}}

\def\BH{\mathcal{B}(H)}

\newcommand{\symm}[1]{{#1}_{\mathrm{sym}}}

\newtheorem{proposition}{Proposition}[section]
\newtheorem{lemma}[proposition]{Lemma}
\newtheorem{theorem}[proposition]{Theorem}

\theoremstyle{definition}

\theoremstyle{remark}
\newtheorem{question}{Question}
\newtheorem{example}[proposition]{Example}

\newtheorem{remark}[proposition]{Remark}
\newtheorem{remarks}[proposition]{Remarks}

\numberwithin{equation}{section}

\begin{document}

\title{Polynomial inequalities for
non-commuting operators}
\author{John E. McCarthy 
\thanks{Partially supported by National Science Foundation Grants DMS 0501079
and DMS 0966845}
\\ Washington University, St. Louis and Trinity College, Dublin \\
\and Richard M. Timoney
\\
Trinity College, Dublin}
\maketitle

\begin{abstract}
We prove an inequality for polynomials applied in a symmetric
way to non-commuting operators.
\end{abstract}

\section{Introduction}

J.~von Neumann \cite{vonNeumann1951}
proved an inequality about the norm of a polynomial
applied to a contraction on a Hilbert space $H$.
Let $\ID$ be the unit disk and $\IT$ the unit circle in $\IC$,
and for any polynomial $p$ let $\| p \|_X $ be the supremum
of the modulus of $p$ on the set $X$.
The
result
is that
\begin{equation}
\label{VNI}
T \in \BH, \|T\| \leq 1 \Rightarrow \|p(T)\| \leq \|p\|_{\CD}.
\end{equation}

For polynomials $p(z) = p(z_1, z_2, \ldots, z_n) = \sum_{|\alpha| \leq N} c_\alpha
z^\alpha $ in $n$ variables we use the standard multi-index notation
(where $\alpha = (\alpha_1, \alpha_2, \ldots, \alpha_n)$ has $0 \leq
\alpha_j \in \IZ$ for $1 \leq j \leq n$,
$|\alpha| = \sum_{j=1}^n \alpha_j$, $z^\alpha = \prod_{j=1}^n z_j
^{\alpha_j}$).
There is an obvious way 
of applying $p$ to an $n$-tuple $T = (T_1, T_2, \ldots, T_n)$
of {\em commuting} operators $T_j \in \BH$ ($1 \leq j \leq n$),
namely
\[
p(T) = p(T_1, T_2, \ldots, T_n)
= \sum_{|\alpha| \leq N} c_\alpha
T^\alpha 
\]
(with $T^\alpha = \prod_{j=1}^n T_j ^{\alpha_j}$ and $T_j^0 = I$).

T. And\^o \cite{Ando1963} proved an extension of von Neumann's inequality to pairs
of commuting contractions.
\begin{theorem}[And\^o]
\label{AndoTheorem}
If $T_1, T_2 \in \BH$, $\max( \|T_1\|, \|T_2\|) \leq 1$, $T_1 T_2 =
T_2 T_1$ and $p(z) = p(z_1, z_2)$ is a polynomial, then
\[
\| p(T_1, T_2) \| \leq \|p\|_{\CD^2}.
\]
\end{theorem}

The purpose of this note is to look for analogues of And\^o's inequality that are 
satisfied by {\em non-commuting} operators. 
For a polynomial $p$ in $n$ variables and  an $n$-tuple of operators $T = (T_1,\dots, T_n)$ 
we define 
$\symm{p}(T)$ 
to be a symmetrized version of $p$ applied to $T$ (we make this precise in Section~\ref{secb}).
We are looking for results of the form:

{\em For all $n$-tuples $T$ of operators in a certain set, there is a set $K_1$ in $\IC^n$ such that
\be
\label{eqa1}
\| \symm{p}(T)  \|
\ \leq \ \| p \|_{K_1} .
\ee
}

and

{\em For all $n$-tuples $T$ of operators in a certain set, there is a set $K_2$ in $\IC^n$ and a constant 
$M$ such that
\be
\label{eqa1a}
\| \symm{p}(T)  \|
\ \leq \ M\, \| p \|_{K_2} .
\ee
}

Our main result is:

 {\bf Theorem~ \ref{thm4k} }
 {\em
There are positive constants $M_n$ and $R_n$
such that, whenever
$T = (T_1, T_2, \ldots, T_n) \in {\BH}^n$ 
satisfies 
\[
\| \sum_{i=1}^n \zeta_i T_i \| \ \leq \ 1 \qquad \forall \zeta_i \in \overline{\ID} ,
\]
and $p$ is a 
polynomial in $n$ variables, then 
\begin{eqnarray}
\label{diam2intro}
\|\symm{p}(T)\|\  & \leq \ & \| p \|_{R_n \CD^n} \\
\label{diam2bintro}
\|\symm{p}(T)\|\  & \leq \ & \ M_n \|p\|_{\CD^n}.
\end{eqnarray}
Moreover, one can choose $R_2 =1.85,\, R_3 = 2.6,\, M_2 = 4.1$ and
$M_3 = 16.6$.
}


\section{Tuples of noncommuting contractions}
\label{secb}

There are several natural ways one might apply a 
polynomial $p(z_1, z_2)$ in two variables to pairs
$T=(T_1, T_2) \in \BH^2$ of
operators.
A simple
case is for polynomials of the form $p(z_1, z_2) = p_1(z_1) +
p_2(z_2)$ where we could naturally consider $p(T_1, T_2)$ to mean 
$p_1(T_1) + p_2(T_2)$.

A recent result of Drury \cite{DruryEJLA2009} is that if $p(z_1, z_2)
= p_1(z_1) +
p_2(z_2)$, $T_1, T_2 \in \BH$ (no longer necessarily commuting),
$\max(\|T_1\|, \|T_2\|) \leq 1$, then
\begin{equation}
\label{eq:drur}
\|p(T_1, T_2)\| \leq \sqrt{2} \|p\|_{\CD^2}.
\end{equation}
Moreover, Drury \cite{DruryEJLA2009} shows that the constant
$\sqrt{2}$ is best possible. 

One way to apply a polynomial $p(z_1, z_2) = \sum_{j,k=0}^n a_j
z_1^j z_2^{k}$ to two noncommuting
operators $T_1$ and $T_2$ is by mapping each monomial $z_1^j z_2^k$ to the
average over all possible products of $j$ number of $T_1$ and $k$
number of
$T_2$, and then extend this map by linearity to all polynomials.
We use the notation $\symm{p}(T_1, T_2)$ and the formula
\[
\symm{p}(T_1, T_2) = \sum_{j,k=0}^n
\frac{a_j}{\binom{j+k}{j}} \sum_{S \in \mathcal{P}(j+k, j)} \prod_{i=1}^{j+k}
T_{2-\chi_S(i)}
\]
where $\mathcal{P}(j+k, j)$ denotes the subsets of $\{1, 2, \ldots,
j+k\}$ of cardinality $j$. The empty product, which arises for $j = k
= 0$, should be taken as the identity operator. The notation
$\prod_{i=1}^{j+k} T_{2-\chi_S(i)}$ is intended to mean the ordered
product
\[
T_{2-\chi_S(1)} T_{2-\chi_S(2)} \cdots T_{2-\chi_S(j+k)},
\]
and $\chi_S(i)$ denotes the indicator function of $S$.

\begin{remarks}
The operation $p \mapsto
\symm{p}(T_1, T_2)$ is not an
algebra homomorphism (from polynomials to operators). It is a linear
operation and does not respect squares in general.

For example, if
$p(z_1, z_2) = z_1^2 + z_2^2$, then 
\[
\symm{p}(T_1, T_2) = T_1^2 +
T_2^2
\]
but for $q(z_1, z_2) = (p(z_1, z_2))^2 = z_1^4 + z_2^4 + 2
z_1^2 z_2^2$ we have 
\[
(\symm{p}(T_1, T_2))^2 = T_1^4 + T_2^4 + T_1^2
T_2^2 + T_2^2 T_1^2 \neq \symm{q}(T_1, T_2)
\]
in general.

Similarly for $p(z_1, z_2) = 2 z_1 z_2$ and 
\[
q(z_1, z_2) = (p(z_1,
z_2))^2 = 4 z_1^2 z_2^2,
\]
$\symm{p}(T_1, T_2) = T_1 T_2 + T_2 T_1$,
\[
(\symm{p}(T_1, T_2))^2 = T_1 T_2 T_1 T_2 + T_1 T_2^2 T_1 + T_2 T_1^2
T_2 + T_2 T_1 T_2 T_1 \neq \symm{q}(T_1, T_2)
\]
in general.

However in the very restricted situation that $p(z_1, z_2) = \alpha +
\beta z_1 + \gamma z_2$ and
$q = p^m$, then we do have
$
\symm{q}(T_1, T_2) = (\symm{p}(T_1, T_2))^m
$.

The symmetrizing idea generalizes in the obvious way to $n > 2$ variables. We will use
the notation $\symm{p}(T)$ for $n$-tuples $T \in \BH^n$ for $n \geq
2$.

\end{remarks}

\section{Example}

The analogue of And\^o's 
inequality for $n \geq 3$  commuting
Hilbert space contractions and polynomials norms on $\ID^n$ is known to
fail (see Varopoulos \cite{VaropoulosJFA1974}, Crabb \& Davie
\cite{CrabbDavieBLMS1975}, Lotto \& Steger
\cite{LottoStegerIIPAMS1994}, Holbrook \cite{Holbrook2001}).

The explicit counterexamples of Kaijser \& Varopoulos
\cite{VaropoulosJFA1974},  and Crabb \& Davie
\cite{CrabbDavieBLMS1975}) have $p(T)$ nilpotent (and so of spectral
radius 0). While the examples of Lotto \& Steger
\cite{LottoStegerIIPAMS1994} and Holbrook \cite{Holbrook2001}) do not
have this property, they are obtained by perturbing examples where
$p(T)$ is nilpotent (and so $p(T)$ has relatively small spectral radius).

It is not known whether there is a constant $C_n $ so that the
multi-variable inequality
\begin{equation}
\label{MultiVNI}
\|p(T)\| = \|p(T_1, T_2, \ldots, T_n)\| \leq C_n \|p\|_{\CD^n}
\end{equation}
holds for all polynomials $p(z)$ in $n$ variables and for all
$n$-tuples $T $
of commuting Hilbert space
contractions. 
%
However, it is well-known
that 
 a spectral radius version of And\^o's inequality is true --- indeed, it holds in any Banach algebra.

\begin{proposition}
\label{spectralradisThm}
If $p$ is a polynomial in $n$ variables and $T = (T_1, T_2, \ldots,
T_n)$ is an $n$-tuple of commuting elements in a Banach algebra, each with norm 
at most one, then 
\begin{equation}
\label{spectralradisuinequality}
\rho(p(T)) = \lim_{m \to \infty} \|(p(T))^m\|^{1/m} \leq \|p\|_{\CD^n}
\end{equation}
\end{proposition}

\begin{proof}
We consider a fixed $n$.
It follows from the Cauchy integral formula, that if $\max_{1 \leq j
\leq n} \|T_j\| \leq r < 1$, then
\begin{equation}
\label{weakVNI}
\|p(T)\| = \|p(T_1, T_2, \ldots, T_n)\| \leq C_r \|p\|_{\CD^n}
\end{equation}
for a constant $C_r$ depending on $r$ (and $n$).

To see this write
\[
p(T) = \frac{1}{(2 \pi i)^n}
\int_{\zeta \in \IT^n} \prod_{j=1}^n p(\zeta) \prod_{j=1}^n
(\zeta_j - T_j)^{-1} \, d \zeta_1 \,  d \zeta_2  \, \ldots \, d\zeta_n
\]
and estimate with the triangle inequality. This shows that $C_r =
(1-r)^{-n}$ will work.

Applying (\ref{weakVNI}) to powers of $p$ and using the spectral
radius formula, we get
\[
\rho(p(T)) \leq \|p\|_{\CD^n},
\]
(provided $\max_{1 \leq j
\leq n} \|T_j\| \leq r < 1$). However, for the general case $\max_{1
\leq j
\leq n} \|T_j\| = 1$, we can apply this to $rT$ to get
\[
\rho(p(T)) = \lim_{r \to 1^-}  \rho(p(rT)) \leq \|p\|_\infty.
\qedhere
\]
\end{proof}

%
%

\begin{example}
\label{ex7}

Let $p(z,w) = (z-w)^2 +2 (z+w) + 1 = z^2 + w^2 - 2 zw  +2 (z+w) + 1$,
\[
T_1 = \begin{pmatrix} \cos(\pi/3) & \sin(\pi/3)\\
\sin(\pi/3) & -\cos(\pi/3) \end{pmatrix} =
\begin{pmatrix} 1/2 & \sqrt{3}/2\\
\sqrt{3}/2 & - 1/2\end{pmatrix},
\]
\[
T_2 = \begin{pmatrix} \cos(\pi/3) & -\sin(\pi/3)\\
-\sin(\pi/3) & -\cos(\pi/3) \end{pmatrix}.
\]

Note that $\|p\|_{{\ID}^2} \geq p(1,-1) = 5$.
To show that
$\|p\|_{{\ID}^2}  \leq 5$, consider the homogeneous polynomial
\[
q(z_1, z_2, z_3) = z_1^2 + z_2^2 + z_3^2 - 2 z_1 z_2 - 2 z_1 z_3 - 2
z_2 z_3 
\]
and observe first that $p(z,w) = q(z, w, -1)$. Moreover
\[
\|p\|_{\ID^2}
= \|p\|_{\IT^2}
=  \|q\|_{\IT^3}
=  \|q\|_{\ID^3},
\]
by homogeneity of $q$ and the maximum principle.
Holbrook \cite[Proposition 2]{Holbrook2001} gives a proof that
$\|q\|_{\ID^3}=5$.

We have 
\begin{eqnarray*}
\symm{p}(T_1, T_2) 
&=& (T_1-T_2)^2 + 2(T_1 + T_2) + I\\
&=& \begin{pmatrix}
0 & \sqrt{3}\\
\sqrt{3} & 0
\end{pmatrix}^2
+ 2 \begin{pmatrix} 1 & 0\\
0 & -1
\end{pmatrix} + I
\\
&=& \begin{pmatrix} 3 & 0\\
0 & 3
\end{pmatrix} + 
\begin{pmatrix} 2 & 0\\
0 & -2
\end{pmatrix}  +
\begin{pmatrix} 1 & 0\\
0 & 1
\end{pmatrix} 
\\
&=&
\begin{pmatrix}6 & 0 \\0 & 2
\end{pmatrix}
\end{eqnarray*}
So $\|\symm{p}(T_1, T_2)\| =6 > 5 = \|p\|_{{\ID}^2}$.
\end{example}

\begin{remark}
The example has hermitian $T_1$ and $T_2$ and a polynomial with real
coefficients and yet $\rho(\symm{p}(T_1, T_2)) > \|p\|_{{\ID}^2}$.
Thus even Proposition~\ref{spectralradisThm} 
does not hold for non-commuting pairs.
\end{remark}

One can show that for the polynomial $p$ of Example~\ref{ex7}, one has the inequality
$$
\| \symm{p} (T_1, T_2) \| \ \leq \ (1 + 4\sqrt{2}) \approx 6.66 
$$
for all contractions $T_1$ and $T_2$.  
This estimate is at least an improvement
over using the sum of the absolute values of the coefficients of $p$, so 
one is led to ask how well can one bound $ \| \symm{p} (T) \|$
for general $p$?


\section{$\left\|\sum \zeta_i T_i \right\| \ \leq \ 1$}
\label{secdiam}

In this section, we shall consider $n$-tuples $T = (T_1,\dots, T_n)$ of operators,
not assumed to be commuting, and we shall make the standing assumption:
\begin{equation}
\label{ineqdiam}
\| \sum_{i=1}^n \zeta_i T_i \| \ \leq \ 1 \qquad \forall \, \zeta_i \, \in \, \CD .
\end{equation}
This will hold, for example, if the condition
\begin{equation}
\label{ineq8}
\sum_{i=1}^n \| T_i \| \ \leq \ 1 
\end{equation}
holds.
We wish to derive bounds on $\| \symm{p}(T)\| $.
We start with the following lemma:

\begin{lemma}
\label{CayleyTransformLemma}
If $S \in \BH$ and $\|S\| < 1$ then
\[
\Re((I +S)(I-S)^{-1}) \geq 0.
\]
\end{lemma}

\begin{proof}
\begin{eqnarray*}
\lefteqn{2 \Re((I +S)(I-S)^{-1})}\\
&=& (I - S^*)^{-1} (I + S^*)
+ (I + S)(I - S)^{-1}\\
&=& (I - S^*)^{-1} \bigl[ (I + S^*)
 (I - S)
+ (I - S^*)(I + S) \bigr]
(I - S)^{-1}\\
&=& 2 (I - S^*)^{-1} [I -  S^*  S)]
(I - S)^{-1}\\
&\geq& 0 .
\end{eqnarray*}
\end{proof}

If $p(z) \, = \, \sum c_\alpha z^\alpha $,
define
\begin{equation}
\label{eqgam}
\Gamma p (z) \ =\  \sum c_\alpha \frac{\alpha !}{|\alpha |!} z^\alpha 
\end{equation}
(as usual, $\alpha!$ means $\alpha_1 ! \cdots \alpha_n !$).
We let $\Lambda$ denote the inverse of $\Gamma$:
\[
\Lambda \sum d_\alpha z^\alpha \ =\ \sum d_\alpha \frac{|\alpha| !}{\alpha !} z^\alpha .
\]

\begin{proposition}
\label{ThmWithL1NormT}
Let $T = (T_1, T_2, \ldots, T_n) \in {\BH}^n$ 
satisfy (\ref{ineqdiam})
and $p(z)$ be a 
polynomial in $n$ variables. Then
\begin{equation}
\label{diam}
\|\symm{p}(T)\|\  \leq \ \|\Gamma p\|_{\CD^n}.
\end{equation}
\end{proposition}

\begin{proof}
We first restrict to the case
\[
\zeta = (\zeta_1, \zeta_2, \ldots, \zeta_n) \in \IT^n \Rightarrow \|\zeta \cdot T\| =
\left\|\sum_{j=1}^n \zeta_j T_j \right\| < 1
\]
and hence
by Lemma~\ref{CayleyTransformLemma}
the operator
\[
(I + \zeta \cdot T) (I - \zeta \cdot T)^{-1}
 =  (I + \zeta \cdot T) \sum_{j=0}^\infty (\zeta \cdot T)^j
= I  + 2 \sum_{j=1}^\infty (\zeta \cdot T)^j
\]
has positive real part
\begin{eqnarray*}
K(\zeta, T) &\ = \ &
\Re \left( (I + \zeta \cdot T) (I - \zeta \cdot T)^{-1} \right)
\\
&=&
I + \sum_{j=1}^\infty (\zeta \cdot T)^j +
\sum_{j=1}^\infty (\bar{\zeta} \cdot T^*)^j
\\
&=&
2 \Re \left[
\sum_{\alpha_1, \dots, \alpha_n = 0}^\infty \frac{|\alpha|!}{\alpha!} \zeta^\alpha
\symm{(z^\alpha)} (T) \right] - I .
\end{eqnarray*}

We can compute that for polynomials $p(z) = p(z_1, z_2, \ldots, z_n)$,
\[
\symm{p}(T) = \int_{\IT^n} \Gamma p(\zeta)  K(\bar\zeta, T)
\, d \sigma(\zeta)
\]
with $d\sigma$ indicating normalised Haar measure on the torus
$\IT^n$ (and $\bar\zeta = (\bar\zeta_1 , \bar\zeta_2, \ldots ,
\bar\zeta_n)$).

As
\[
K(\bar\zeta, T) \, d \sigma(\zeta)
\]
is a positive operator valued measure on $\IT^n$,
we then have
a positive unital linear map $C(\IT^n) \to \BH$ given by
$
f \mapsto  \int_{\IT^n} f(\zeta) K(\bar\zeta, T)
\, d \sigma(\zeta).
$
As this map is then of norm 1, we can conclude
\[
\|\symm{p}(T)\| \leq \|\Gamma p\|_{\CD^n}.
\]

For the remaining case $
\sup_{\zeta \in \IT^n} \|\zeta \cdot T\| =1$,
we have
\[
\|\symm{p}(T)\| = \lim_{r \to 1^-} \| \symm{p}(rT)\| \leq \| \Gamma p\|_{\CD^n}.
\qedhere
\]
\end{proof}

\begin{remark}
The technique of the above proof is derived from methods of
\cite{McCarthyPutinar2005}.
\end{remark}

Now we want to estimate $\| \Gamma p \|_{\CD^N}$.

\begin{proposition}
\label{mspec}
For each $n \geq 2$ there is a constant $M_n$ so that
$$
\| \Gamma p \|_{\CD^n} \ \leq \ M_n \, \| p \|_{\CD^n} .
$$
Moreover,
\begin{eqnarray*}
M_2 &\ \leq \ & 4.07 \\
M_3 &\leq &
16.6
\end{eqnarray*}
\end{proposition}
\begin{proof}
Define
\begin{equation}
\label{eq4h}
J(\eta) \ = \ \sum_{\alpha_1 =0,\dots,\alpha_n=0}^\infty \frac{\alpha!}{|\alpha|!} \eta^\alpha .
\end{equation}
Then
\begin{equation}
\label{eq4a}
\Gamma p (z) \ = \ \int_{\IT^n} p(\zeta) [ J(z_1 \bar \zeta_1, \dots, z_n \bar \zeta_n)] d\sigma(\zeta) .
\end{equation}
To use (\ref{eq4a}), we break $J$ into two parts --- the sum $J_0$ where the
minimum of the $\alpha_i$ is $0$, and  the remaining terms $J_1$.
\[
J_1(\eta) \ = \ \sum_{\alpha_1 =1,\dots,\alpha_n=1}^\infty \frac{\alpha!}{|\alpha|!} \eta^\alpha .
\]

Case: $n=2$. Here,
\begin{equation}
\label{eq4b}
\int_{\IT^2} p(\zeta) J_0 ( z_1 \bar \zeta_1, z_2 \bar \zeta_2) d\sigma(\zeta) \ = \
p(z_1,0) + p(0,z_2) - p(0,0) .
\end{equation}
So the norm of the left-hand side of (\ref{eq4b}) is dominated by $3\| p\|_{\CD^2}$.

For $J_1$, we will use the estimate
\[
\left| \int_{\IT^2} p(\zeta) J_1 ( z_1 \bar \zeta_1, z_2 \bar \zeta_2) d\sigma(\zeta) \right|
\ \leq \ \| p \|_{\infty} \| J_1 \|_{L^1} \ \leq \
\| p \|_{\infty} \| J_1 \|_{L^2} .
\]
We have
\begin{eqnarray*}
\| J_1 \|^2_{L^2} & \ = \ &
\sum_{\alpha_1, \alpha_2 =1}^\infty \left( \frac{\alpha_1 ! \alpha_2 !}{(\alpha_1 + \alpha_2)!} \right)^2
\\
&\ = \ &
\sum_{\alpha_1 =1}^\infty \frac{1}{(\alpha_1 +1)^2} \, + \,
\sum_{\alpha_2 =2}^\infty \frac{1}{(\alpha_2 +1)^2} \, + \,
\sum_{\alpha_1, \alpha_2 =2}^\infty \left( \frac{\alpha_1 ! \alpha_2 !}{(\alpha_1 + \alpha_2)!} \right)^2
\\
&\ \leq \ &
\left( \frac{\pi^2}{3} - \frac{9}{4} \right) \, + \,
\sum_{k=4}^\infty (k-3) \left( \frac{2}{k(k-1)} \right)^2 \\
&\ \leq \ & (1.069)^2 .
\end{eqnarray*}
(In the penultimate line, we let $k = \alpha_1 + \alpha_2$; there are $k-3$ terms with this
sum, and the largest they can be is when either $\alpha_1$ or $\alpha_2$ is $2$.)
Adding the two estimates, we get $M_2 \leq 4.07$.

Case: $n=3$. Again, we estimate the contributions of $J_0$ and $J_1$ separately.
We have
\begin{eqnarray*}
\lefteqn{\int p (\zeta) J_0 (z_1 \bar \zeta_1,z_2 \bar \zeta_2,z_3 \bar \zeta_3) d\sigma (\zeta) }
\\
& &= \
\Gamma p(0,z_2, z_3) + [\Gamma p(z_1,0,z_3) - p(0,0,z_3)] \\
&&\ + \, [\Gamma p(z_1,z_2,0) - p(z_1,0,0) - p(0,z_2,0) + p(0,0,0)]
\end{eqnarray*}
where we have had to subtract some terms to avoid double-counting.
Thus the contribution of $J_0$ is at most
$3 M_2 + 4$.

To calculate the contribution of $J_1$, we make the following estimate on
$\| J_1 \|_{L^2}$, which is valid for all $n \geq 3$:

We want to bound
\begin{equation}
\label{eq4c}
\sum_{\alpha_1 =1, \dots, \alpha_n =1}^\infty \left( \frac{\alpha !}{|\alpha|!} \right)^2
\end{equation}
Let $k = |\alpha|$ in (\ref{eq4c}).
Note first that the number of terms for each $k$ is the number of ways of writing $k$ as a sum
of $n$ distinct positive integers (order matters),
and this is exactly $\binom{k-1}{n-1}$.
Moreover, as each $\alpha_i$ is at least $1$, we have
\[
\frac{\alpha !}{|\alpha|!}
\ \leq \
\frac{1}{k(k-1) \cdots (k -n +2)} .
\]
Therefore (\ref{eq4c}) is bounded by
\begin{eqnarray*}
\label{eq4d}
\lefteqn{
\sum_{k=n}^\infty \binom{k-1}{n-1} \left(
\frac{1}{k(k-1) \cdots (k -n +2)} \right)^2
} \\
&& = \
\sum_{k=n}^\infty
\frac{k-n+1}{(n-1)! k}
\frac{1}{k(k-1) \cdots (k -n +2)} .
\end{eqnarray*}

The terms on the right-hand side of (\ref{eq4d}) decay like $1/k^{n-1}$,
so the series converges for all $n \geq 3$.
When $n=3$, the series is
\[
\sum_{k=3}^\infty
\frac{k-2}{2k^2 (k-1)} \ \leq \ (0.381)^2.
\]
Therefore $M_3 \leq 3M_2 + 4.381 \, < 16.59$.

We now proceed by induction on $n$.
The contribution from $J_0$ is dominated by applying $\Gamma$ to the restriction of
$p$ to the slices with one or more coordinates equal to $0$, and these are bounded by
the inductive hypothesis. The contribution from $J_1$ is bounded by (\ref{eq4c}).

\end{proof}

We have proved that the polydisk is an $M$-spectral set for $T$; we can make the
constant one by enlarging the domain.

\begin{proposition}
\label{propspset}
There is a constant $R_n$ so that
\begin{equation}
\label{eq4g}
\| \Gamma p \|_{\CD^n} \ \leq \
\| p \|_{R_n \CD^n} .
\end{equation}
Moreover,
\begin{eqnarray*}
R_2 &\ \leq \ & 1.85 \\
R_3 &\leq &
2.6
\end{eqnarray*}
\end{proposition}
\begin{proof}
Let $L(\eta) = 2 \Re J(\eta) -1$.
Adding terms that are not conjugate analytic powers of $\zeta$ inside the bracket in (\ref{eq4a}) will not change the value of the integral,
so, writing $z \bar \zeta$ for the $n$-tuple $(z_1 \bar \zeta_1, \dots, z_n \bar \zeta_n)$, we get
\begin{equation}
\label{eq4f}
\Gamma p (z) \ = \ \int_{\IT^n} p(\zeta) [ L(z \bar \zeta )] d\sigma(\zeta) .
\end{equation}
As $L$ is real and has integral $1$, if we can choose $r_n$ so that 
if $| z_i | \leq r_n$ for each $i$ then $L(z \bar \zeta)$ is non-negative for all $\zeta$,
then its $L^1$ norm would equal its integral, and so we would get from (\ref{eq4f}) that
\[
| \Gamma p (z) | \ \leq \ \| p \|_{\CD^n} .
\]
Letting $R_n = 1/r_n$ gives (\ref{eq4g}).
As the series (\ref{eq4h})
converges absolutely for all $\eta \, \in \, \ID^n$, and
$L(0) =1$, the existence of some $r_n$ now follows by continuity.

Let us turn now to obtaining  quantitative estimates.

Case: $n=2$.
Adding terms to $J$ that are not analytic will not affect the integral
(\ref{eq4f}), so let us consider 
\[
L'(\eta) 
\ = \ 
\Re\left[ \frac{1 + \eta_1}{1-\eta_1} \right] \, \cdot \,
\Re\left[ \frac{1 + \eta_2}{1-\eta_2} \right]
\ - \ 
\sum_{\alpha_1 =1, \alpha_2 =1}^\infty
( 1 - \frac{\alpha!}{|\alpha|!} )
(\eta_1^{\alpha_1} - \bar \eta_1^{\alpha_1} )
(\eta_2^{\alpha_2} - \bar \eta_2^{\alpha_2} )
.
\]
Then $L'$ has integral $1$ and (\ref{eq4f}) is unchanged if $L$ is replaced by $L'$.
So we wish to find the largest $r$ so that
$L'$ is positive on $r \ID^2$.

It can be checked numerically that $r=0.5406$ works, so 
the best $R_2$ is smaller than the reciprocal of 0.5406, which is less than 1.85.

Case: $n=3$. (By hand).

Our strategy will be to simply estimate each non-constant term in $L$ by a function of $r$, add them up, and
see how small $r$ must be for all these terms to be less than $1$.

First, let us estimate the terms from $J_0 -1$, i.e. those terms with either one or two of the $\alpha_i$'s equal
to $0$.
We can write
\[
J_0 (\eta) -1 \ = \ [ \eta_1 + \sum_{\alpha_2 =1, \alpha_3 =1}^\infty \frac{\alpha_2 ! \alpha_3 !}{(\alpha_2 + \alpha_3)!} 
\eta_2^{\alpha_2} \eta_3^{\alpha_3} ] + \dots ,
\]
where the $\dots$ mean two more terms with the indices $(1,2,3)$ permuted.
Therefore if each $|\eta_i | \leq r$, we have
\begin{eqnarray*}
2 \Re [ J_0 (\eta) -1 ]
&\ \leq \ &
6r \, + \, 6  \sum_{\alpha_2 =1, \alpha_3 =1}^\infty \frac{\alpha_2 ! \alpha_3 !}{(\alpha_2 + \alpha_3)!} r^{\alpha_2 + \alpha_3} 
\\
&\leq&
6r \, + \, 6 (r^2/2 + 2r^3/3) +  6 \sum_{k=4}^\infty r^k (k-1) \frac{1}{k(k-1)} 
\\
&=&
6r \, + \, 6 (r^2/2 + 2r^3/3) +  6 [- \log(1-r) -r - r^2/2 -r^3/3 ] .
\end{eqnarray*}

The contribution to $L$ from $J_1$, where all the indices are at least $1$, is at most
\begin{eqnarray*}
2 \Re [ J_1 (\eta)  ]
&\ \leq \ &
2  \sum_{\alpha_1 =1, \alpha_2 =1, \alpha_3 =1}^\infty \frac{\alpha ! }{|\alpha|!} r^{|\alpha|} 
\\
&\leq&
2 \sum_{k=3}^\infty r^k \binom{k-1}{2} \frac{1}{k(k-1)} 
\\
&=&
 \frac{r^3}{1-r} + 2 [ \log (1-r) + r + r^2/2] .
 \end{eqnarray*}
Adding the two terms together, we get
\[
L(\eta) \ \geq \ 
1 - [ 2r + r^2 + 2r^3 + \frac{r^3}{1-r} - \log(1-r) ] ,
\]
and this is positive if $r \leq .152$.
So letting  $R_3$ be the reciprocal of this root, which is less than 6.6, will work.  

Case: $n=3$. (Computer-aided)

As in the  case $n=2$, we consider the kernel
\begin{eqnarray*}
L'(\eta) 
&\ =\ & 
\Re\left[ \frac{1 + \eta_1}{1-\eta_1} \right] \, \cdot \,
\Re\left[ \frac{1 + \eta_2}{1-\eta_2} \right] \, \cdot \,
\Re\left[ \frac{1 + \eta_3}{1-\eta_3} \right]
\\ 
&&\ - \ 
\sum_{\alpha_1 =1, \alpha_2 =1, \alpha_3 =0}^\infty
( 1 - \frac{\alpha!}{|\alpha|!} )
(\eta_1^{\alpha_1} - \bar \eta_1^{\alpha_1} )
(\eta_2^{\alpha_2} - \bar \eta_2^{\alpha_2} )
(\eta_3^{\alpha_3} + \bar \eta_3^{\alpha_3} )
.
\end{eqnarray*}
(Note that there is a plus in the last factor to keep $L'$ real.)
Again, a computer search can find $r$ so that $L'$ is positive on $r \ID^3$, and
$r=.39$ works, so $R_3 < 2.6$.

\end{proof}

Combining Propositions~\ref{ThmWithL1NormT}, \ref{mspec}
and \ref{propspset}, we get the main result of this section.

\begin{theorem}
\label{thm4k}
There are positive constants $M_n$ and $R_n$ such that
whenever
$T = (T_1, T_2, \ldots, T_n) \in {\BH}^n$ 
satisfies (\ref{ineqdiam})
and $p(z)$ is a 
polynomial in $n$ variables, then 
\begin{eqnarray}
\label{diam2}
\|\symm{p}(T)\|\  & \leq \ & \| p \|_{R_n \CD^n} \\
\label{diam2b}
\|\symm{p}(T)\|\  & \leq \ & \ M_n \|p\|_{\CD^n}.
\end{eqnarray}
Moreover, one can choose $R_2 =1.85,\, R_3 = 2.6,\, M_2 = 4.1$ and
$M_3 = 16.6$.
\end{theorem}

\begin{remark}
Another way to estimate $\| \symm{p} (T) \|$,
under the assumption (\ref{ineq8}), would be to crash through with absolute values.
Let
$\Delta_n = \{ z \in \IC^n : \sum_{j=1}^n |z_j| \leq 1\}$
and let $r_n$ denote the Bohr radius of $\Delta_n$, {\em i.e.}
the largest $r$ such that whenever
$p(z) \, = \, \sum c_\alpha z^\alpha$ has modulus one on $\Delta_n$,
then
$q(z) \, = \, \sum |c_\alpha| z^\alpha$ has modulus bounded by one
on $r \Delta_n$.
One then has the estimate that, under the hypothesis 
(\ref{ineq8}),
and writing $C_n = 1/r_n$,
\begin{equation}
\label{diam3}
\|\symm{p}(T)\| \ \leq\  \| q \|_{\Delta_n} \ \leq \ 
\|p\|_{C_n \, \Delta_n}.
\end{equation}
It was shown by 
L.~Aizenberg \cite[Thm. 9]{aiz00} that 
$$
\frac{1}{3 e^{1/3}}
\ < \
r_n   
\ \leq \
\frac{1}{3}.
$$
So 
the estimate in (\ref{diam2}) 
for pairs satisfying (\ref{ineq8}) 
does not follow from (\ref{diam3}).
\end{remark}

\section{$n$-tuples of contractions}
\label{sec:con}

In an attempt to use the above technique for tuples $T \in \BH^n$ such
that $\max_{1 \leq j \leq n} \|T_j\| \leq 1$, we consider restricting
$\zeta$ to belong to
$\Delta_n$, and we replace $\sigma$ by some probability measure $\mu$ supported on $\Delta_n$.

Suppose we can find some function $q$ such that
\begin{equation}
\label{eq45}
\Lambda_\mu(q)(z) \ :=\  \int_{\Delta_n} q(\zeta) \Re \frac{ 1 +
\bar{\zeta} \cdot z}{ 1 - 
\bar{\zeta} \cdot z} \, d\mu(\zeta)
\end{equation}
equals $p(z)$.
We do not actually need $q$ to be a polynomial; having an absolutely convergent power
series on $\Delta_n$ (in $\zeta$ and $\bar \zeta$) is enough.


\begin{lemma}
\label{leme1}
With notation as above, assume $\Lambda_\mu(q) = p$
and that  $T \in \BH^n$ is an $n$-tuple of contractions. Then
\[
\|\symm{(  p )}(T) \| \leq \|q\|_{{\rm suppt}(\mu)} \leq \  \sup \{
|q(z)| : z \in \Delta_n\}.
\]
\end{lemma}

\begin{proof}
We assume first that $\max_{1 \leq j \leq n} \|T_j\| < 1$ and use the
notation $K(\zeta, T)$ from the proof of Proposition~\ref{ThmWithL1NormT}
(which is permissible as $\|\zeta \cdot T\| < 1$ for $\zeta \in
\Delta_n$).
We have
\[
\symm{( \Lambda_\mu q )}(T) = \int_{\Delta_n} q(\zeta) K \bar\zeta,
T) \, d \sigma(\zeta)
\]
and hence the inequality $\|\symm{(  p )}(T) \| \leq
\|q\|_{{\rm suppt} (\mu) }$ follows as in the previous proof.

If $\max_{1 \leq j \leq n} \|T_j\| =1$, we deduce the result from
$\|\symm{(  p )}(r T) \| \leq
\|q\|_{\Delta_n}$ for $0 < r < 1$.
\end{proof}

\begin{remark}
For an arbitrary measure $\mu$, there might be no $q$ such that
$\Lambda_\mu (q) = p$. If $\mu$ is chosen to be circularly symmetric, though, one gets
\begin{equation}
\label{eq73}
\Lambda_\mu (z^\alpha) \ = \
\left[ \frac{|\alpha|!}{\alpha_1 ! \dots \alpha_n !} \int |\zeta^\alpha|^2 d\mu(\zeta) \right] z^\alpha.
\end{equation}
As long as none of the moments on the right of (\ref{eq73}) vanish, inverting $\Lambda_\mu$ is
now straightforward.
\end{remark}

To make use of the lemma to bound $\symm{p}(T)$ we need to find a way
to choose another polynomial $q$ and a 
$\mu$ on $\Delta_n$  so that $p = \Lambda_\mu q$ and
$\|q\|_{\Delta_n}$ is small.
We do not know a good way to do this.

\begin{question} What is the smallest constant $R_n$ such that, for
every $n$-tuple $T$ of contractions and every polynomial $p$, one has
\be
\label{que1}
\|\symm{p}(T)\|\  \leq \| p \|_{ R_n \CD^n} ?
\ee
\end{question}
We do not know if one can choose $R_n$ smaller than the reciprocal of the
Bohr radius of the polydisk, even when $n=2$.

\begin{question} Is there a  constant $M_n$ such that, for
every $n$-tuple $T$ of contractions and every polynomial $p$, one has
\be
\label{que2}
\|\symm{p}(T)\|\  \leq M_n \, \| p \|_{  \CD^n} ?
\ee
\end{question}

\bibliography{vni}
\bibliographystyle{rt}

\end{document}